\documentclass[11pt]{article}
\usepackage{geometry}

\newtheorem{theorem}{Theorem}

\newtheorem{conjecture}[theorem]{Conjecture}

\newtheorem{proposition}[theorem]{Proposition}

\newenvironment{proof}[1][Proof]{\noindent\textbf{#1.} }{\ \rule{0.5em}{0.5em}}
\input{tcilatex}

\begin{document}

\title{On the X-rays of Permutations}
\author{Cecilia Bebeacua \\
Imperial College \and Toufik Mansour \\
Univ. of Haifa \and Alexander Postnikov \\
MIT \and Simone Severini \\
Univ. of York}
\maketitle

\begin{abstract}
The X-ray of a permutation is defined as the sequence of antidiagonal sums
in the associated permutation matrix. X-rays of permutation are interesting
in the context of Discrete Tomography since many types of integral matrices
can be written as linear combinations of permutation matrices. This paper is
an invitation to the study of X-rays of permutations from a combinatorial
point of view. We present connections between these objects and
nondecreasing differences of permutations, zero-sum arrays, decomposable
permutations, score sequences of tournaments, queens' problems and rooks'
problems.
\end{abstract}

\section{Introduction}

Let $\mathcal{S}_{n}$ be the set of all permutations of $[n]=\{1,2,\ldots,n%
\} $ and let $P_{\pi}$ be the permutation matrix corresponding to $\pi \in%
\mathcal{S}_{n}$. For $k=2,\ldots,2n$, the $(k-1)$-th \emph{antidiagonal sum}
of $P_{\pi}$ is $x_{k-1}(\pi)={\textstyle\sum_{i+j=k}}[P_{\pi}]_{i,j}$. The
sequence of nonnegative integers $x(\pi)=x_{1}(\pi )x_{2}(\pi)\ldots
x_{2n-1}(\pi)$ is called the (\emph{antidiagonal}) \emph{X-ray} of $\pi$.
The \emph{diagonal X-ray} of $\pi$, denoted by $x_{d}(\pi)$, is similarly
defined. Note that $x(\pi)=x(\pi^{-1})$, for every $\pi\in\mathcal{S}_{n}$.
The sequence $x(\pi)$ may be also seen as a word over the alphabet $[n]$. As
an example, the following table contains the X-rays of all permutations in $%
\mathcal{S}_{3}$:\vspace{10pt}

{\small 
\begin{tabular}{||c|c||c|c||c|c||c|c||c|c||}
\hline
$\pi$ & $x(\pi)$ & $\pi$ & $x(\pi)$ & $\pi$ & $x(\pi)$ & $\pi$ & $x(\pi)$ & $%
\pi$ & $x(\pi)$ \\ \hline
$123$ & $10101$ & $231,312$ & $01110$ & $132$ & $10020$ & $213$ & $02001$ & $%
321$ & $00300$ \\ \hline
\end{tabular}%
}\vspace{10pt}

Although X-rays of permutations are interesting object on their own, among
the reasons why they are of general interest in Discrete Tomography \cite{k}
is that many types of integral matrices can be written as linear
combinations of permutation matrices (for example, binary matrices with
equal row-sums and column-sums, like the adjacency matrices of Cayley
graphs). Deciding whether for a given word $w=w_{1}\ldots w_{2n-1}$ there
exists $\pi\in\mathcal{S}_{n}$ such that $w=x(\pi)$ is an NP-complete
problem \cite{b} (see also \cite{g}). The complexity is polynomial if the
permutation matrix is promised to be wrapped around a cylinder \cite{d}. It
is necessary to keep into account that permutations are not generally
specified by their X-rays: just consider the permutation $\pi=73142865$ and $%
\sigma=72413865$; we have $x(\pi )=x(\sigma)=000110200002100$, $%
x_{d}(\pi)=x_{d}(\sigma)=00021111100010$ and $\pi\neq\sigma^{-1}$. This
hints that an issue concerning X-rays is to quantify how much information
about $\pi$ is contained in $x(\pi)$. In this paper we present some
connections between X-rays of permutations and a variety of combinatorial
objects. From a practical perspective, this may be useful in isolating and
approaching special cases of the above problem.

The remainder of the paper is organized as follows. In Section 2 we consider
the problem of counting X-rays. We prove a bijection between X-rays and
nondecreasing differences of permutations. We define the \emph{degeneracy}
of an X-ray $x(\pi)$ as the number of permutations $\sigma$ such that $%
x(\pi)=x(\sigma)$, and we characterize the X-rays with the maximum
degeneracy. We prove a bijection between X-ray of length $4k+1$ having
maximum degeneracy and zero-sum arrays. In Section 3 we consider the notion
of simple permutations. This notion seems to provide a good framework to
study the degeneracy of X-rays, but the relation between simple permutations
and X-rays with small degeneracy remains unclear. Section 4 is devoted to
binary X-rays, that is X-rays whose entries are only zeros and ones. We
characterize the X-rays of circulant permutation matrices of odd order.
Moreover, we present a relation between binary X-rays, the $n$-queens
problem (see, \emph{e.g.}, \cite{v}), the score sequences of tournaments on $%
n$ vertices (see \cite[Sequence A000571]{s}), and extremal Skolem sequences,
see~\cite[Conjecture 2.2]{N}.

A number of conjectures and open problems will be explicitly formulated or
will simply stand out from the context. We use the standard notation for
integers sequences from the OEIS \cite{s}.

\section{Counting X-rays}

We begin by addressing the following natural question:\ what is the number
of different X-rays of permutations in $\mathcal{S}_{n}$? Although we are
unable to find a generating function for the sequence, we show a bijection
between X-rays and nondecreasing differences of permutations. The \emph{%
difference} of permutations $\pi,\sigma\in\mathcal{S}_{n}$ is the integers
sequence $\pi-\sigma=(w_{1},w_{2},\ldots,w_{n})$, where $w_{1}=\pi_{1}-%
\sigma_{1},w_{2}=\pi_{2}-\sigma_{2},\ldots,w_{n}=\pi_{n}-\sigma_{n}$. For
example, if $\pi=1234$ and $\sigma=2413$, we have $e-2413=(-1,-2,2,1)$. Let $%
x_{n}$ be the numbers of different X-rays of permutations in $\mathcal{S}%
_{n} $. Let $d_{n}$ be the number of nondecreasing differences of
permutations in $\mathcal{S}_{n}$. The number $d_{n}$ equals the number of
different differences $e-\sigma$ with entries rearranged in the
nondecreasing order. In other words, $d_{n}$ equals the number of different
multisets of the form $M(\sigma
)=\{1-\sigma_{1},2-\sigma_{2},\ldots,n-\sigma_{n}\}$, with entries
rearranged in the nondecreasing order. The entries of $x(\pi)$ are then the
entries of the vector $e_{1-\sigma_{1}}+e_{2-\sigma_{2}}+\ldots+e_{n-%
\sigma_{n}}$, where $e_{i}$ is the $i$-th coordinate vector of length $2n-1$%
. For example, for $\pi=3124$ we have $x(3124)=0101200$ and $%
e_{1-3}+e_{2-1}+e_{3-2}+e_{4-4}=(0,1,0,0,0,0,0)+(0,0,0,0,1,0,0)+(0,0,0,0,1,0,0)+(0,0,0,1,0,0,0)=(0,1,0,1,2,0,0) 
$. On the basis of this reasoning we can state the following result.

\begin{proposition}
The number $x_{n}$ of different X-rays of permutations in $\mathcal{S}_{n}$
is equal to the number $d_{n}$ \textrm{(}see \textrm{\cite[Sequence A019589]%
{s})} of nondecreasing differences of permutations in $\mathcal{S}_{n}$.
\end{proposition}

Let us define and denote the \emph{degeneracy} of an X-ray $x(\pi)$ by 
\[
\delta(x(\pi))=|\{\sigma:x(\sigma)=x(\pi)\}|. 
\]
If $x(\pi)$ is such that $\delta(x(\pi))\geq\delta(x(\sigma))$ for all $%
\sigma\in S_{n}$, we write $x_{\max}^{n}=x(\pi)$ and we say that $x(\pi)$
has \emph{maximum degeneracy}. The following table contains $x_{n}$, $%
x_{\max}^{n}$ and $\delta(x_{\max}^{n})$ for $n=1,\ldots,8$.\vspace{10pt}

\begin{tabular}{||llll||llll||}
\hline
$n$ & $x_{n}$ & $x_{\max}^{n}$ & $\delta(x_{\max}^{n})$ & $n$ & $x_{n}$ & $%
x_{\max}^{n}$ & $\delta(x_{\max}^{n})$ \\ \hline
$1$ & \multicolumn{1}{|l}{$1$} & \multicolumn{1}{|l}{$1$} & 
\multicolumn{1}{|l||}{$1$} & $5$ & \multicolumn{1}{|l}{$59$} & 
\multicolumn{1}{|l}{$001111100$} & \multicolumn{1}{|l||}{$6$} \\ \hline
$2$ & \multicolumn{1}{|l}{$2$} & \multicolumn{1}{|l}{$020,101$} & 
\multicolumn{1}{|l||}{$1$} & $6$ & \multicolumn{1}{|l}{$246$} & 
\multicolumn{1}{|l}{$00011211000$} & \multicolumn{1}{|l||}{$12$} \\ \hline
$3$ & \multicolumn{1}{|l}{$5$} & \multicolumn{1}{|l}{$01110$} & 
\multicolumn{1}{|l||}{$2$} & $7$ & \multicolumn{1}{|l}{$1105$} & 
\multicolumn{1}{|l}{$0001111111000$} & \multicolumn{1}{|l||}{$28$} \\ \hline
$4$ & \multicolumn{1}{|l}{$16$} & \multicolumn{1}{|l}{$0012100$} & 
\multicolumn{1}{|l||}{$3$} & $8$ & \multicolumn{1}{|l}{$5270$} & 
\multicolumn{1}{|l}{$000011121110000$} & \multicolumn{1}{|l||}{$76$} \\ 
\hline
\end{tabular}%
\vspace{15pt}

It is not difficult to characterize the X-rays with maximum degeneracy. One
can verify by induction that for $n$ even,%
\[
x_{\max}^{n}=00\ldots011\ldots121\ldots110\ldots00, 
\]
with $n/2$ left-zeros and right-zeros, and $n/2-1$ ones; for $n$ odd, 
\[
x_{\max}^{n}=00..011\ldots110..00, 
\]
with $(n-1)/2$ left-zeros and right-zeros, and $n$ ones. Notice that if $%
x(\pi)=x_{\max}^{n}$ (for $n$ odd) then $P_{\pi}$ can be seen as an
hexagonal lattice with all sides of length $\left( n+1\right) /2$. In each
cell of the lattice there is $0$ or $1$, and $1$ is in exactly $n$ cells;
the column-sums are $1$ and the diagonal and anti-diagonal sums are $0$.
This observation describes a bijection between permutations of odd order
whose X-ray is $x_{\max}^{n}$ and zero-sum arrays. An $(m,2n+1)$-\emph{%
zero-sum} array is an $m\times(2n+1)$ matrix whose $m$ rows are permutations
of the $2n+1$ integers $-n,-n+1,\ldots,n$ and in which the sum of each
column is zero \cite{bp}. The matrix 
\[
\left[ 
\begin{array}{rrr}
-1 & 0 & 1 \\ 
0 & 1 & -1 \\ 
1 & -1 & 0%
\end{array}
\right] 
\]
is an example of $(3,3)$-zero-sum array. Thus we have the next result.

\begin{proposition}
The number $\delta(x_{\max}^{n})$ for $n$ odd is equal to the number of $%
(3,2n+1)$-zero-sum arrays \textrm{(}see \textrm{\cite[Sequence A002047]{s})}.
\end{proposition}

Before concluding the section, it may be interesting to notice that if we
sum entry-wise the X-rays of all permutations in $\mathcal{S}_{n}$ we obtain
the following sequence of $2n-1$ terms: 
\[
(n-1)!,2(n-1)!,\ldots,(n-1)(n-1)!,n!,(n-1)(n-1)!,\ldots,2(n-1)!,(n-1)!. 
\]
The meaning of the terms of this sequence is clear.

\section{Simple permutations and X-rays}

In the previous section we have considered the X-rays with maximum
degeneracy. What can we say about X-rays with degeneracy $1$? If $%
\delta(x(\pi))=1$ then $\pi$ is an involution (in such a case $%
P_{\pi}=P_{\pi}^{-1}$) but the converse if not necessarily true. In fact
consider the involution $\pi=1267534$. One can verify that $%
x(\pi)=x(\sigma)=x(\rho)=1010000212000$, for $\rho=1275634$ and $%
\sigma=1267453$. In a first approach to the problem, it seems useful to
study what kind of operations can be done \textquotedblleft
inside\textquotedblright\ a permutation matrix $P_{\pi}$ in order to obtain
another permutation, say $P_{\sigma}$, such that $x(\pi)=x(\sigma)$ and $%
P_{\pi}\neq P_{\sigma}^{-1}$. A intuitively good framework for this task is
provided by the notion of \emph{block permutation}. A \emph{segment} and a 
\emph{range} of a permutation are a set of consecutive positions and a set
of consecutive values. For example, in the permutation $34512$, the segment
formed by the positions $2,3,4$ is occupied by the values $4,5,1$; the
elements $1,2,3$ form a range. A \emph{block} is a segment whose values form
a range. Every permutation has singleton blocks together with the block $%
12\ldots n$. A permutation is called \emph{simple} if these are the only
blocks \cite{a}. A permutation is said to be a \emph{block permutation} if
it is not simple. Note that if $\pi$ is simple then it is $\pi^{-1}$. Let $%
S=(\pi_{1}\in\mathcal{S}_{n_{1}},\ldots ,\pi_{k}\in\mathcal{S}_{n_{k}})$ be
an ordered set and let $\pi\in\mathcal{S}_{k}$. We assume that in $S$ there
exists $1\leq i\leq k$ such that $n_{i}>1$. We denote by $P(\pi,S)$ the $%
(n_{1}+\cdots +n_{k})$-dimensional permutation matrix which is partitioned
in $k^{2}$ blocks, $B_{1,1},\ldots ,B_{k,k}$, such that $B_{i,j}=P_{\pi_{i}}$
if $\pi(i)=j$ and $B_{i,j}=0$, otherwise. We denote by $\pi\lbrack\pi_{1},%
\ldots ,\pi_{k}]$ (or equivalently by $(\pi)[S]$) the permutation
corresponding to $P(\pi,S)$. For example, let $S=(231,21,312)$ and $\pi=231$%
. Then 
\[
P(231,S)=\left[ 
\begin{array}{ccc}
0 & P_{231} & 0 \\ 
0 & 0 & P_{21} \\ 
P_{312} & 0 & 0%
\end{array}
\right] 
\]
and $231[S]=231[231,21,312]=56487312$. The matrix $P(231,S)$ can be modified
leaving the X-ray of $(\pi)[S]$ invariant:%
\[
P(231,(312,21,312))=\left[ 
\begin{array}{ccc}
0 & P_{312}=P_{231}^{T} & 0 \\ 
0 & 0 & P_{21} \\ 
P_{312} & 0 & 0%
\end{array}
\right]. 
\]
It is clear that $56487312$ is a block permutation. Let $\pi$ be a simple
permutation then possibly $\delta(x(\pi))>1$. In fact, the permutation $%
\pi=531642$ is simple, but $\delta(x(\pi))=6$, since $x(\pi
)=00111011100=x(526134)=x(461253)$, plus the respective inverses. The
permutations $526134$ and $461253$ are decomposable. This means that there
possibly exists a decomposable permutation $\sigma$ such that $x(\sigma
)=x(\pi)$, even if $\pi$ is simple. There relation between simple
permutations and X-rays of small degeneracy is not clear. Intuitively, a
simple permutation allows less \textquotedblleft freedom of
movement\textquotedblright\ than a block permutation. It is also intuitive
that we have low degeneracy when the nonzero entries of the X-ray are
\textquotedblleft distributed widely\textquotedblright\ among the $2n-1$
coordinates. The following result is easily proved.

\begin{proposition}
Let $\sigma=\pi\lbrack S]=\pi\lbrack\pi_{1},\ldots,\pi_{k}]$ be a block
permutation. Then $\delta(x(\pi))>1$ if one of the following two conditions
is satisfied:

(1) If $\pi\neq12\ldots n$ then there is at least one $\pi_{i}\in S$ which
is not an involution;

(2) If $\pi=12\ldots n$ then there are at least two $\pi_{i},\pi_{i}\in S$
which are not involution.
\end{proposition}

\begin{proof}
(1) Let $\pi\neq12\ldots n$ be any permutation. Take $\pi_{i}^{-1}$ for some 
$\pi_{i}\in S$. Let $\rho=\pi\lbrack\pi_{1},\ldots,\pi_{i}^{-1},\ldots,%
\pi_{k}]$. Since $\sigma$ is a block permutation, $x(\sigma)=x(\rho)$.
However, if $\pi_{i}\neq\pi_{i}^{-1}$ then $\sigma\neq\rho$ and $%
\sigma^{-1}\neq\sigma$. It follows that $x(\sigma)$ does not specify $\sigma$%
. (2) Let $\pi=12\ldots n$. Let all elements of $S$ be involutions except $%
\pi_{i}$. Take $\pi_{i}^{-1}$. Let $\rho=\pi\lbrack\pi
_{1},\ldots,\pi_{i}^{-1},\ldots,\pi_{k}]$. Again, $x(\sigma)=x(\rho)$, but
this time $\rho=\sigma^{-1}$. Then $x(\sigma)$ possibly specifies $\sigma$.
If, for distinct $i,j$, there are $\pi_{i},\pi_{j}\in S$ such that $%
\pi_{i}\neq\pi _{i}^{-1}$ and $\pi_{j}\neq\pi_{j}^{-1}$ then 
\[
x(\sigma^{\prime})=x(\pi
\lbrack\pi_{1},\ldots,\pi_{i}^{-1},\ldots,\pi_{j}^{-1},\ldots,\pi_{k}])=x(%
\sigma),
\]
but $x(\sigma)$ does not specify $\sigma$, given that $\rho\neq\sigma^{-1}$.
\end{proof}

This is however not a sufficient condition for having $\delta(x(\pi))>1$.
Permutations with equal X-rays are said to be in the same \emph{degeneracy
class}. The table below contain the number of permutations in $\mathcal{S}%
_{n}$ which are in each degeneracy class, and the number of different
degeneracy classes with the same cardinality, for $n=2,\ldots ,8$. These
numbers provide a partition on $n!$. We denote by $C(n)$ the total number of
degeneracy classes. We write $a(b)$, where $a$ is the number of permutations
in the degeneracy class and $b$ the number of degeneracy classes of the same
cardinality:\vspace{10pt}

\begin{tabular}{||l||}
\hline
$C(2)=1$: 1(2) \\ \hline
$C(3)=2$: 1(4),2(1) \\ \hline
$C(4)=3$: 1(9),2(6),3(1) \\ \hline
$C(5)=5$:\ 1(20),2(26),3(6),4(6),6(1) \\ \hline
$C(6)=10$: 1(49),2(100),3(19),4(43),5(1),6(19),7(2),8(11),9(1),2(1) \\ \hline
$C(7)=20$: 1(114),2(345),3(60),4(229),5(18),6(118),7(11),8(98),10(29) \\ 
\hline
11(2),12(33),14(13),16(14),18(6),20(4),21(1),22(2),26(1),28(1). \\ \hline
\end{tabular}%
\ .\vspace{10pt}

We conjecture that if $\delta(x(\pi))=1$ then $x(\pi)$ does not have more
than $2$ adjacent nonzero coordinates. However the converse is not true if $%
\pi \in\mathcal{S}_{n}$ for $n\geq8$: for $\pi=17543628$ and $%
\sigma=16547328 $, we have $x(\pi)=x(\sigma)=100000320010001$, but there are
no more than $2$ adjacent coordinates.

\section{Binary X-rays}

In general, it does not seem to be an easy task to characterize X-rays. A
special case is given by X-rays associated with circulant permutation
matrices, for which is available an exact characterization. An X-ray $x(\pi)$
is said to be \emph{binary} if $x_{i}(\pi)\in\{0,1\}$ for every $1\leq
i\leq2n-1$. The set all permutations in $\mathcal{S}_{n}$ with binary X-ray
is denoted by $\mathcal{B}_{n}$. Counting binary X-rays means solving a
modified version of the $n$-queens problem (see, \emph{e.g.}, \cite{v}) in
which two queens do not attack each other if they are in the same
NorthWest-SouthEst diagonal. The permutations with binary X-rays associated
to circulant matrices are characterized in a straightforward way. Let $C_{n}$
be the permutation matrix associated with the permutation $c_{n}=23\ldots n1$%
, that is the \emph{basic circulant permutation matrix}. The matrices in the
set $\mathcal{C}_{n}=\{C_{n}^{0},C_{n},C_{n}^{2},\ldots,C_{n}^{n-1}\}$ ($%
C_{n}^{0}$ is the identity matrix) are called the \emph{circulant
permutation matrices}. The matrix $C_{n}^{k}$ is associated to $c_{n}^{k}$.
Observe that $x(\pi)$ can be seen as a binary number, since $%
x_{i}(\pi)\in\{0,1\}$ for every $i$. Let 
\[
d_{j}(\pi)=2^{2n-1-j}\cdot x_{j}(\pi),\quad j=1,2,\ldots,2n-1, 
\]
and $d(\pi)=\sum_{i=1}^{2n-1}d_{i}(\pi)$, that is the decimal expansion of $%
x(\pi)$. The table below lists the X-rays of $\mathcal{C}_{3},\mathcal{C}%
_{5} $ and $\mathcal{C}_{7}$, and their decimal expansions:\vspace{10pt} 
\[
\begin{tabular}{||l|l|l||l|l|l||}
\hline
$\pi$ & $x(\pi)$ & $d(\pi)$ & $\pi$ & $x(\pi)$ & $d(\pi)$ \\ \hline
$123$ & $10101$ & $21$ & $12345$ & $101010101$ & $341$ \\ \hline
$231$ & $01110$ & $14$ & $23451$ & $010111010$ & $186$ \\ \hline
&  &  & $34512$ & $001111100$ & $124$ \\ \hline
\end{tabular}%
\ . 
\]

For $\pi=c_{n}^{k}$, one can verify that 
\[
\begin{array}{ll}
d(\pi) & =\frac{1}{6}2^{\frac{3}{2}n+\frac{1}{2}+k}-\frac{1}{6}2^{\frac {1}{2%
}n+\frac{1}{2}+k}+\frac{1}{3}2^{\frac{3}{2}n+\frac{1}{2}-k}-\frac{1}{3}2^{%
\frac{1}{2}n+\frac{1}{2}-k} \\ 
& =a(k)\left( 2^{n}-1\right) \left( 2^{n}-1\right) 2^{\frac{1}{2}n-k+\frac{1%
}{2}},%
\end{array}%
\]
where $a(k)=(2^{2k-1}+1)/3$ (A007583).

In the attempt to count binary X-rays, we are able to establish a bijection
between these objects and score sequences of tournaments. A \emph{tournament}
is a loopless digraph such that for every two distinct vertices $i$ and $j$
either $\left( i,j\right) $ or $\left(j,i\right)$ is an arc \cite{bk}. The 
\emph{score sequence} of an tournament on $n$ vertices is the vector of
length $n$ whose entries are the out-degrees of the vertices of the
tournament rearranged in nondecreasing order.

\begin{proposition}
Let $b_{n}$ be the number of binary X-rays of permutations in $S_n$ and let $%
s_{n}$ be the number of different score sequences of tournaments on $n$
vertices \textrm{(}see \textrm{\cite[Sequence A000571]{s})}. Then $b_n\leq
s_n$.
\end{proposition}

\begin{proof}
The number $s_{n}$ equals the number of integers lattice points $%
(p_{0},\dots,p_{n})$ in the polytope $P_{n}$ given by the inequalities $%
p_{0}=p_{n}=0$, $2p_{i}-p_{i+1}-p_{i-1}\leq1$ and $p_{i}\geq0$, for $%
i=1,\dots,n-1$, see \cite{bk}. Let $x_{1},\dots,x_{n}$ be the coordinates
related to $p_{1},\dots,p_{n}$ by $p_{i}=x_{1}+\dots+x_{i}-i^{2}$, for $%
i=1,\dots,n$. We can rewrite the inequalities defining the polytope $P_{n}$
in these coordinates as follows: $x_{1}+\dots+x_{i}\geq i^{2}$, $x_{i+1}\geq
x_{i}+1$ and $x_{1}+\dots+x_{n}=n^{2}$. For a permutation $w\in\mathcal{S}%
_{n}$ with a binary X-ray, let $l_{i}=l_{i}(w)$ be the position of the $i$%
-th `$1$' in its X-ray. In other words, the sequence $(l_{1},\dots,l_{n})$
is the increasing rearrangement of the sequence $(w_{1},w_{2}+1,w_{3}+2,%
\dots,w_{n}+n-1)$. Then the numbers $l_{1},\dots ,l_{n}$ satisfy the
inequalities defining the polytope $P_{n}$ (in the $x$-coordinates). Indeed, 
$l_{1}+\dots+l_{n}=w_{1}+(w_{2}+1)+\cdots +(w_{n}+n-1)=n^{2}$; $l_{i+1}\geq
l_{i}+1$; and the minimal possible value of $l_{1}+\cdots+l_{i}$ is $%
(1+0)+(2+1)+\cdots+(i+(i-1))=i^{2}$. This finishes the proof. In order, to
prove that $b_{n}=s_{n}$ it is enough to show that, for any integer point $%
(x_{1},\dots,x_{n})$ satisfying the above inequalities, we can find a
permutation $w\in\mathcal{S}_{n}$ with $x_{i}=l_{i}(w)$.
\end{proof}

\begin{conjecture}
\label{conjbin} All binary X-rays of permutations in $\mathcal{S}_{n}$ are
in a bijective correspondence with integer lattice points $%
(x_{1},\dots,x_{n})$ of the polytope given by the inequalities 
\[
\begin{array}{l}
x_{1}+\cdots+x_{i}\geq i^{2},\qquad i=1,\dots,n; \\ 
x_{1}+\cdots+x_{n}=n^{2}, \\ 
x_{i+1}-x_{i}\geq1,\qquad i=1,\dots,n-1.%
\end{array}
\]
For a permutation $w\in S_{n}$, the corresponding sequence $(x_{1},\dots
,x_{n})$ is defined as the increasing rearrangement of the sequence $%
(w_{1},w_{2}+1,w_{3}+2,\dots,w_{n}+n-1)$.
\end{conjecture}

Again, it is clear that X-rays injectively map into the integer points of
the above polytope. One needs to show that there will be no gaps in the
image. Also, it can be shown that the above conjecture is equivalent to
Conjecture 2.2 from \cite{N} concerning extremal Skolem sequences. The
conjecture turns out to be false, when not restricted to binary X-rays.

We conjecture also that the number of different X-rays of permutations in 
\emph{$\mathcal{S}_{n}$} whose possible entries are $0$ and $2$ is equal to
the number of score sequences in tournament with $n$ players, when 3 points
are awarded in each game (see \cite[Sequence A047729]{s}).

\section{Palindromic X-rays}

What can we say about X-rays with special symmetries? The \emph{reverse} of $%
x(\pi)$, denoted by $\overline{x}(\pi)$, is the mirror image of $x(\pi)$. If 
$x(\pi)=\overline{x}(\pi)$ then $\pi$ is said to be \emph{palindromic}. The 
\emph{reverse} of $\pi$, denoted by $\overline{\pi}$, is mirror image of $%
\pi $. For example, if $\pi=25143$ then $\overline{\pi}=34152$. The
permutation matrix $P_{\overline{\pi}}$ is obtained by writing the rows of $%
P_{\pi}$ in reverse order. In general $\overline{x}(\pi)\neq x(\overline{\pi}%
)$. In fact, for $\pi=25143$, we have $x(\pi)=0011001200$, $\overline{x}%
(\pi)=0021001100$ and $x(\overline{\pi})=0020011010$. We denote by $|M$ and $%
\underline{M}$ the matrices obtained by writing the columns and the rows of
a matrix $M$ in reverse order, respectively.

\begin{proposition}
Let $l_{n}$ be the number of permutations in $\mathcal{S}_{n}$ with
palindromic X-rays and let $i_{n}$ be the number of involutions in $\mathcal{%
S}_{n}$ \textrm{(}see \textrm{\cite[Sequence A000085]{s})}. Then, in
general, $l_{n}>i_{n}$.
\end{proposition}

\begin{proof}
Recall that a permutation $\pi$ is an \emph{involution} if $\pi=\pi^{-1}$.
Since $P_{\pi}=P_{\pi}^{T}$, it is clear that the diagonal X-ray of an
involution $\pi$ is palindromic. The X-ray of $\sigma$ such that $P_{\sigma
}=|P_{\pi}$ is then also palindromic. This shows that $l_{n}\geq i_{n}$.
Now, consider a permutation matrix of the form%
\[
P_{\sigma}=\left[ 
\begin{array}{cc}
P_{\rho} & 0 \\ 
0 & P_{\rho}^{T}%
\end{array}
\right], 
\]
for some permutation $\rho$ which is not an involution. Then $P_{\rho}\neq
P_{\rho}^{T}$, $P_{\sigma}\neq P_{\sigma}^{T}$ and $\sigma$ is not an
involution, but the diagonal X-ray of $\sigma$ is palindromic. The X-ray of $%
\pi$ such that $P_{\pi}=|P_{\sigma}$ is then also palindromic. This proves
the proposition.
\end{proof}

The next contains the values of $l_{n}$ for small $n$: 
\[
\begin{tabular}{||l|l||l|l||l|l||l|l||}
\hline
$n$ & $l_{n}$ & $n$ & $l_{n}$ & $n$ & $l_{n}$ & $n$ & $l_{n}$ \\ \hline
$2$ & $2$ & $4$ & $12$ & $6$ & $128$ & $8$ & $2110$ \\ \hline
$3$ & $4$ & $5$ & $32$ & $7$ & $436$ & $9$ & $8814$ \\ \hline
\end{tabular}
\ . 
\]

\begin{proposition}
Let $l_{n,A=D}$ be the number of permutations in $\mathcal{S}_{n}$ with:

(1) equal diagonal and antidiagonal X-rays;

(2) palindromic X-rays.\newline
Let $r_{n}$ be the number of permutations in $\mathcal{S}_{n}$ invariant
under the operation of first reversing and then taking the inverse \textrm{(}%
see \textrm{\cite[Sequnce A097296]{s})}. Then, in general, $l_{n,A=D}>r_{n}$.
\end{proposition}

\begin{proof}
We first construct the permutations which are invariant under the operation
of first reversing and then taking the inverse. Let $\pi\in\mathcal{S}_{n}$
where $n=2k$. We look at $P_{\pi}$ as partitioned in $4$ blocks: 
\[
P_{\pi}=\left[ 
\begin{array}{cc}
A & B \\ 
C & D%
\end{array}
\right] . 
\]
If 
\[
P_{\pi}=(\underline{P_{\pi}})^{T}=\left[ 
\begin{array}{cc}
A & B \\ 
C & D%
\end{array}
\right] ^{T}=\left[ 
\begin{array}{cc}
\underline{B} & \underline{A} \\ 
\underline{D} & \underline{C}%
\end{array}
\right] ^{T}=\left[ 
\begin{array}{cc}
(\underline{B})^{T} & (\underline{D})^{T} \\ 
(\underline{A})^{T} & (\underline{C})^{T}%
\end{array}
\right] 
\]
then $A=(\underline{B})^{T},B=(\underline{D})^{T},C=(\underline{A})^{T}$ and 
$D=(\underline{C})^{T}$. This implies the X-ray of $P_{\pi}$ being
palindromic and, moreover, the diagonal and antidiagonal X-rays being equal.
Note that we can construct $P_{\pi}$ only if $n\equiv0(\func{mod}4)$, and in
this case $r_{n}\neq0$. However, fixed $n\equiv0(\func{mod}4)$, we have $%
r_{n}=r_{n+1}$, since the permutation matrix%
\[
P_{\sigma}=\left[ 
\begin{array}{ccc}
A &  & \mathbf{0} \\ 
& 1 &  \\ 
\mathbf{0} &  & D%
\end{array}
\right] +\left[ 
\begin{array}{ccc}
\mathbf{0} &  & B \\ 
& 1 &  \\ 
C &  & \mathbf{0}%
\end{array}
\right] 
\]
can be always constructed from $P_{\pi}$. (Permutation matrices like $P_{\pi}
$ and $P_{\sigma}$ provide the solutions of the \textquotedblleft
rotationally invariant\textquotedblright\ $n$-rooks problem. This points out
that A097296 and A037224 are indeed the same sequence.) Now, the proposition
is easily proved by observing that, for $\rho=369274185$, $P_{\rho}$ is not
of the form of $P_{\sigma}$. A direct calculation shows that $r_{9}=12$ and $%
l_{9,A=D}=20$.
\end{proof}


\begin{thebibliography}{99}
\bibitem{a} Albert, M. H., M. D. Atkinson, and M. Klazar, The enumeration of
simple permutations, \emph{Journal of Integer Sequences} \textbf{6} (2003),
Article 03.4.4, 18 pp.

\bibitem{bp} Bennett, B. T. and R. B. Potts, Arrays and brooks, \emph{J.
Austral. Math. Soc. }\textbf{7} (1967), 23--31.

\bibitem{b} Brunetti, S., A. Del Lungo, P. Gritzmann, and S. de Vries, On
the reconstruction of permutation and partition matrices under tomographic
constraints, \emph{preprint}.

\bibitem{d} Del Lungo, A., Reconstructing permutation matrices from diagonal
sums. Selected papers in honour of Maurice Nivat. \emph{Theoret. Comput. Sci.%
} \textbf{281:1-2} (2002), 235--249.

\bibitem{g} Gardner, R. J., P. Gritzmann, and D. Prangenberg, On the
computational complexity of reconstructing lattice sets from their X-rays, 
\emph{Discrete Math.} \textbf{202} (1999), 45-71.

\bibitem{k} Herman, G. T. and A. Kuba (Eds.),\emph{\ Discrete tomography:
Foundations, Algorithms and Applications,} Birkh\"{a}user Boston, 1999.

\bibitem{N} Nordh, G., Generalization of Skolem Sequences, Master's Thesis,
LITH-MAT-EX-2003/05, Department of Mathematics, Link\"opings universitet,
2003.

\bibitem{bk} Reid, K. B., Tournaments, \emph{The Handbook of Graph Theory,}
J. L. Gross and J. Yellen (editors), CRC Press, Boca Raton, Fl., 156-184,
2004.

\bibitem{v} Ruskey, F., 
\texttt{http://www.theory.csc.uvic.ca/\symbol{126}cos/inf/misc/Queen.html}%
.

\bibitem{s} Sloane, N. J. A., (2005), The On-Line Encyclopedia of Integer
Sequences, 
\texttt{http://www.research.att.com/\symbol{126}njas/sequences/}%
.
\end{thebibliography}
\end{document}